\documentclass[reqno]{amsart}
\usepackage[foot]{amsaddr}

\usepackage{type1cm}        

\usepackage{graphicx}       
\usepackage{multicol} 
\usepackage[bottom]{footmisc}

\usepackage{amsmath}
\usepackage{amsthm}
\usepackage{amssymb}
\usepackage{makecell}
\usepackage{multirow}
\usepackage{tcolorbox}
\tcbuselibrary{skins}

\usepackage{hyperref}

\newtheorem{observation}{Observation}
\newtheorem{lemma}{Lemma}

\newcommand{\bx}{\mathbf{x}}
\newcommand{\ba}{\mathbf{a}}
\newcommand{\bg}{\mathbf{g}}
\newcommand{\bX}{\mathbf{X}}
\newcommand{\bhX}{{\widehat\bX}}
\newcommand{\hw}{{\widehat w}}
\newcommand{\ub}{\text{\rm ub}}
\newcommand{\lb}{\text{\rm lb}}
\newcommand{\oR}{\mathbb{R}}
\newcommand{\oN}{\mathbb{N}}
\newcommand{\oS}{\mathbb {S}}
\newcommand{\fmin}{f_{\text{\rm min}}}
\newcommand{\fminX}{f_{\text{\rm min},\bX}}

\newcommand{\gminX}{g_{\text{\rm min},\bX}}
\newcommand{\gminhX}{g_{\text{\rm min}, \bhX}}

\newcommand{\MP}{{\mathcal P}}
\newcommand{\MQ}{{\mathcal Q}}
\newcommand{\MT}{{\mathcal T}}
\newcommand{\MS}{{\mathcal S}}
\newcommand{\MM}{{\mathcal M}}
\newcommand{\MB}{{\mathcal B}}
\newcommand{\MC}{\mathcal{C}}

\newcommand{\mc}{\text{\rm mc}}
\newcommand{\kernelop}{\mathbf{K}}
\newcommand{\kernel}{\mathcal{K}}
\newcommand{\by}{\mathbf{y}}

\newcommand{\Se}{\mathrm{Se}}
\newcommand{\bh}{\mathbf{h}}
\newcommand{\bY}{\mathbf{Y}}

\newcommand{\preo}[1]{\mathcal{T}(#1)}
\newcommand{\qmod}[1]{\mathcal{Q}(#1)}

\newcommand{\mainset}{\mathbf{X}}
\newcommand{\x}{\bx}
\newcommand{\ball}[1]{B^{#1}}
\newcommand{\simplex}[1]{\Delta^{#1}}

\newcommand{\Tr}{\text{\rm Tr}}
\newcommand{\sfT}{{\sf T}}
\newcommand{\ignore}[1]{}

\newenvironment{trailer}[1]{\begin{tcolorbox}[enhanced, left=0pt, title=\textbf{{#1}}, colback=white, colbacktitle=gray!20, coltitle=black, sharp corners, frame hidden, borderline south = {1pt}{1pt}{gray!20}]}{\end{tcolorbox} \medskip\noindent}

\calclayout

\begin{document}

\title[An Overview of Convergence Rates for Sum of Squares Hierarchies]{An Overview of Convergence Rates for Sum of Squares Hierarchies in Polynomial Optimization}

\author{Monique Laurent$^1$}
\address{$^1$CWI Amsterdam and Tilburg University}
\author{Lucas Slot$^2$}
\address{\quad $^2$ETH Zurich}
\email{m.laurent@cwi.nl, lucas.slot@inf.ethz.ch}

\begin{abstract} In this survey we consider polynomial optimization problems, asking to minimize a polynomial function over a compact semialgebraic set, defined by polynomial inequalities. This models a great variety of (in general, nonlinear nonconvex) optimization problems. Various hierarchies of (lower and upper) bounds have been introduced, having the remarkable property that they converge asymptotically to the global minimum. These bounds  exploit algebraic representations of positive polynomials in terms of sums of squares and can be computed using semidefinite optimization.  
Our focus lies in the performance analysis of these hierarchies of bounds, namely, in  how far the bounds are from the global minimum as the degrees of the sums of squares they involve tend to infinity. We present the main known state-of-the-art results and offer a gentle introductory overview over the various  techniques that have been recently developed to establish them, stemming from the theory of orthogonal polynomials, approximation theory, Fourier analysis, and more.
\end{abstract}

\maketitle

\section{Introduction}

This survey offers a gentle introduction and overview over the {\em design and performance analysis of approximation hierarchies for  polynomial optimization}. In this section  we will first introduce polynomial optimization, and its use for modeling hard optimization problems arising within a broad range of fields and application domains. After that we will explain how to design hierarchies of approximations using sums of squares of polynomials as a tractable surrogate for polynomial nonnegativity, and how sums of squares of polynomials can be modeled using semidefinite optimization.

\subsection{Polynomial Optimization}
Throughout we use the following notation. We let $\oR[\bx]=\oR[x_1,\ldots,x_n]$ denote the ring of $n$-variate real polynomials in the variables $\bx=(x_1,\ldots,x_n)$.  For a multi-index $\alpha\in \oN^n$, we let $x^\alpha=x_1^{\alpha_1}\cdots x_n^{\alpha_n}$ denote the associated monomial, whose degree is $|\alpha|:=\|\alpha\|_1=\alpha_1+\ldots + \alpha_n$. We let $\oN^n_d$ denote the set of sequences $\alpha\in \oN^n$ with $|\alpha|\le d$.
Any  polynomial $f\in \oR[\bx]$ can be expressed in the monomial basis as $f=\sum_\alpha f_\alpha x^\alpha$ for some coefficients $f_\alpha\in \oR$. Here, the sum is finite and the largest value of $|\alpha|$ for which $f_\alpha\ne 0$ is the {\em degree} of $f$, denoted $\deg(f)$. For an integer $d\ge 0$, we let $\oR[\bx]_d$ denote the set of polynomials with degree at most $d$. It is a vector space, with 
$[x]_d=\{x^\alpha: \alpha\in \oN^n_d\}$
as its standard monomial basis. As we will see later, it is sometimes convenient to use other polynomial bases, that are orthogonal with respect to some inner product induced by a selected measure on   the space $\oR^n$.

The general setting of polynomial optimization is as follows. We are given polynomials $f,g_1,\ldots,g_m\in \oR[\bx]$ and the goal is to find the minimum value that $f$  takes on the feasible region defined by the polynomials $g_j$. In other words, the goal is to  solve
the following optimization problem.
\begin{trailer}{The polynomial optimization problem}
\begin{equation}\label{eqfmin}
\begin{split}
\fmin& = \min\{f(\bx): g_1(\bx)\ge 0,\ldots, g_m(\bx)\ge 0\}\\
& =\min\{f(\bx): \bx \in \bX\},
\end{split}
\end{equation}
after setting 
\begin{equation}\label{eqX}
\bX=\{\bx\in \oR^n: g_1(\bx)\ge 0,\ldots, g_m(\bx)\ge 0\}.
\end{equation}
We also sometimes use the notation $\fminX$ for the minimum value of (\ref{eqfmin}) and the notation $\Se(\bg)$ for the set in \eqref{eqX}. 
\end{trailer}
Such a set $\bX$ is  known as a basic closed semialgebraic set. Typical instances include the unit sphere $\oS^{n-1}$, the standard simplex $\Delta^n$, the unit ball $B^n$, the box $[-1,1]^n$ or $[0,1]^n$, and the binary cube $\{\pm 1\}^n$ or $\{0,1\}^n$, where
\begin{equation*}\begin{split}
\oS^{n-1}=\{\bx\in \oR^n: \|\bx\|=1\},\ B^n=\{\bx\in\oR^n: \|\bx\|\le 1\},\\
 \Delta^n=\Big\{\bx\in \oR^n: x_i\ge 0 \ (i\in [n]), \sum_{i=1}^n x_i=1\Big\},
 \end{split}
\end{equation*}
and $\|\bx\|^2={\sum_{i=1}^nx_i^2}$ denotes the squared Euclidean norm. Throughout we assume $\bX$ is compact. So, problem (\ref{eqfmin}) always has a global minimizer.

The polynomial optimization problem (\ref{eqfmin}) contains linear programming as a (very) special case, when all polynomials $f,g_1,\ldots,g_m$ are linear (i.e., of degree at most~$1$). In this case, efficient optimization algorithms exist (see, e.g.,  Roos et al. \cite{RoosTerlakyVial2005}). But, it also models a much broader range of problems, which are generally nonlinear and nonconvex. This includes well-known NP-hard problems, already when restricting to seemingly simple feasible regions $\bX$ in (\ref{eqX}), such as the simplex, the box, or the sphere as we illustrate on some examples in the next section.

Polynomial optimization has received growing research interest in the past decades, when it was realized that algebraic and geometric properties of polynomials could be exploited to design dedicated methods, able to capture the {\em global minimum}, in contrast to general nonlinear optimization methods where one can often only gain information about local minima. In a nutshell, this research direction builds on combining real algebraic geometry results (about sums of squares of polynomials) and functional analytical results (about moments of measures) with semidefinite optimization. It roots in foundational works, in particular, by Shor~\cite{Shor1987}, Nesterov~\cite{Nesterov2000}, Lasserre~\cite{Lasserre2001}, Parrilo~\cite{Parrilo2000, Parrilo2003}. The field has substantially grown\footnote{This is witnessed, e.g., by the fact that it has  received its own Mathematics Subject Classification number 90C23.} and  has a broad literature. We  mention some books and overviews that can serve as introduction to the topic and give further references to many applications and additional aspects that are not mentioned in the present paper; in particular, by  Lasserre~\cite{Lasserre2009,Lasserre2015}, Laurent~\cite{Laurent2009}, Blekherman, Parrilo, Thomas~\cite{BPTbook}, Henrion, Korda, Lasserre~\cite{HKLbook}, Magron, Wang~\cite{Magron:sparsebook}, Nie~\cite{Niebook}.
  
The present paper will focus on the performance analysis of various hierarchies of bounds that have been introduced for the polynomial optimization problem \eqref{eqfmin}, based on using tailored sums of squares representations for positive polynomials.  We will recall the definition of these hierarchies of (upper and lower) bounds and discuss the main state-of-the-art results that have been shown in recent years about their quantitative convergence properties (these are not covered in the literature mentioned above). We focus on offering a gentle overview of the main techniques that are needed to prove these quantitative results.
    
\subsection{Examples and Applications}

We begin with mentioning a few instances of the polynomial optimization problem~\eqref{eqfmin} that capture well-known hard combinatorial optimization problems. 

Consider  a graph $G=(V,E)$, where $V=[n]=\{1,\ldots,n\}$ is the set of vertices  and the pairs in $E\subseteq V\times V$ correspond to the edges of $G$. A set $I\subseteq V$ is said to be independent (or stable) if it contains no edge and a fundamental combinatorial problem is determining the largest cardinality of an independent set, denoted $\alpha(G)$, a well-known NP-hard problem (see Garey, Johnson \cite{GJ1979}). 
Interestingly, this problem admits several equivalent reformulations as  instances of polynomial optimization over the boolean cube, the box, the simplex, the sphere, respectively:
\begin{align*}
\alpha(G) & =\max\Big\{\sum_{i\in V} x_i: x_ix_j=0 \text{ for } \{i,j\}\in E, x_i^2=x_i \text{ for } i\in V\Big\},\\
\alpha(G)& =\max\Big\{ \sum_{i\in V}x_i-\sum_{\{i,j\}\in E} x_ix_j: x\in [0,1]^n\Big\},\\
{1\over \alpha(G)}& =\min\Big\{x^\sfT (I+A_G)x: x\in \oR^n_+, \sum_{i\in V}x_i =1\Big\},\\
{1\over \alpha(G)}& =\min\Big\{(x^{\circ 2})^\sfT (I+A_G)x^{\circ 2}:
x\in \oR^n, \sum_{i\in V}x_i^2=1\Big\},
\end{align*}
where we set $x^{\circ 2}=(x_1^2,\ldots,x_n^2)^\sfT$ and $A_G$ is the adjacency matrix of  graph $G$. An additional reformulation  is in terms of linear optimization over copositive matrices:
\begin{align*}
\alpha(G)=\min\{\lambda: \lambda(I+A_G)-J \text{ is copositive}\},
\end{align*}
where $J$ denotes the all-ones matrix. We refer to de Klerk, Pasechnik \cite{dKP2002} for details about the above  formulations.

\ignore{
\begin{itemize}
\item[(1)] \ Optimization over the boolean cube:
$$\alpha(G)=\max\Big\{\sum_{i\in V} x_i: x_ix_j=0 \text{ for } \{i,j\}\in E, x_i^2=x_i \text{ for } i\in V\Big\}.$$
\item[(2)] \ Optimization over the box:
$$\alpha(G)=\max\Big\{ \sum_{i\in V}x_i-\sum_{\{i,j\}\in E} x_ix_j: x\in [0,1]^n\Big\}.$$
\item[(3)] \ Optimization over the simplex, with $A_G$ denoting the adjacency matrix of $G$:
$${1\over \alpha(G)}=\min\Big\{x^\sfT (I+A_G)x: x\in \oR^n_+, \sum_{i\in V}x_i =1\Big\},$$
\item[(4)]
\ Optimization over the sphere:
$${1\over \alpha(G)}=\min\Big\{\sum_{i\in V}x_i^4 +2\sum_{\{i,j\}\in E} x_i^2x_j^2: \sum_{i\in V}x_i^2=1\Big\}.$$
\item[(5)] \ Optimization over the cone of copositive matrices, with $J$ denoting the all-ones matrix:
$$\alpha(G)=\min\{\lambda: \lambda(I+A_G)-J \text{ is copositive}\}.$$
\end{itemize}
We refer to de Klerk, Pasechnik \cite{dKP2002} for details about the formulations in (3)-(5).
}

Another fundamental combinatorial problem is the maximum cut problem, asking for the largest cardinality of a cut in a graph $G=(V,E)$, denoted as $\mc(G)$, that can be formulated as any of the following polynomial optimization problems:
\begin{align*}
\mc(G)=\max\Big\{\sum_{\{i,j\}\in E} {1\over 2}(1-x_ix_j): x\in \{\pm 1\}^n\Big\},\\
\mc(G)=\max\Big \{ {1\over 4} x^\sfT L_G x: x\in [-1,1]^n\Big\},
\end{align*}
where $L_G=D-A_G$ is the Laplacian matrix, with $D$ the diagonal matrix having the degrees of the vertices as diagonal entries. 

A fundamental question in analysis  that can be formulated in terms of polynomial optimization is testing whether an $n$-variate  polynomial $f$ is convex. Indeed, $f$ is convex if and only if its Hessian matrix 
$H_f(\bx)=\Big({\partial^2f\over \partial x_i\partial x_j}(\bx)\Big)_{i,j=1}^n $ is positive semidefinite at any $\bx\in \oR^n$ or, equivalently, if the $2n$-variate polynomial $F(\bx,\by)=\by^\sfT H_f(\bx)\by$ is nonnegative over $\oR^n\times \oR^n$ (i.e., its global minimum is 0). It has been shown by Ahmadi et al. \cite{Ahmadi2013} that testing whether a quartic polynomial is convex is an NP-hard problem. This implies hardness of testing global nonnegativity of a quartic polynomial.

So, the above examples show that computing the global minimum of a polynomial over simple regions such as the simplex, the sphere, the box, the full space, or the boolean box, are computationally hard already when restricting to small degree (at most 4).  Another notable application which can be cast as instance of polynomial optimization problem involving only quadratic polynomials, is the optimal power flow problem in energy (see, e.g., Zohrizadeh et al. \cite{Josz}). 
Polynomial optimization has a broad modeling power and permits to capture problems from various areas, such as probablity, mathematical finance, control, game theory, for which  we refer to the exposition by  Lasserre \cite{Lasserre2009}.

\subsection{Nonnegative Polynomials and Sums of Squares}

A polynomial $f\in \oR[\bx]_{2d}$ is called a {\em sum of squares} (abbreviated as {\em sos})
if it can be written as a sum of squares of other polynomials, i.e., $f=\sum_{j=1}^k q_j^2$ for some $q_j\in \oR[\bx]$ and $k\ge 1$. Then, each $q_j$ has degree at most $d$ and can be assumed to be homogeneous if $f$ is homogeneous.
We let $\Sigma$ denote the set of sums of squares of polynomials and set $\Sigma_{2d}=\Sigma\cap \oR[\bx]_{2d}$.  We may use the notation $\Sigma[\bx]$ to stress which  variables are used.

Let $\MP$ denote the set of polynomials $f$ that are globally nonnegative (i.e., $f(\bx)\ge 0$ for all $\bx\in\oR^n$),   and $\MP(\bX)$ the set of polynomials that are nonnegative over a given set $\bX\subseteq \oR^n$. 

Clearly, every sum of squares of polynomials is globally nonnegative, i.e., the inclusion
$\Sigma\subseteq \MP$ holds. As is well-known, this inclusion is in general strict. 
By a result of Hilbert \cite{Hilbert1888} we know that any globally nonnegative $n$-variate polynomial with degree $2d$ is a sum of squares of polynomials only in the following three  cases: $n=1$ (univariate), $d=1$ (quadratic), and the exceptional case $(n=2,d=2)$ (i.e., quartic in two variables).

\begin{figure}
\hspace*{-0.5cm}
\parbox{7cm}{
\vspace*{-5cm}
\includegraphics[width=7cm]{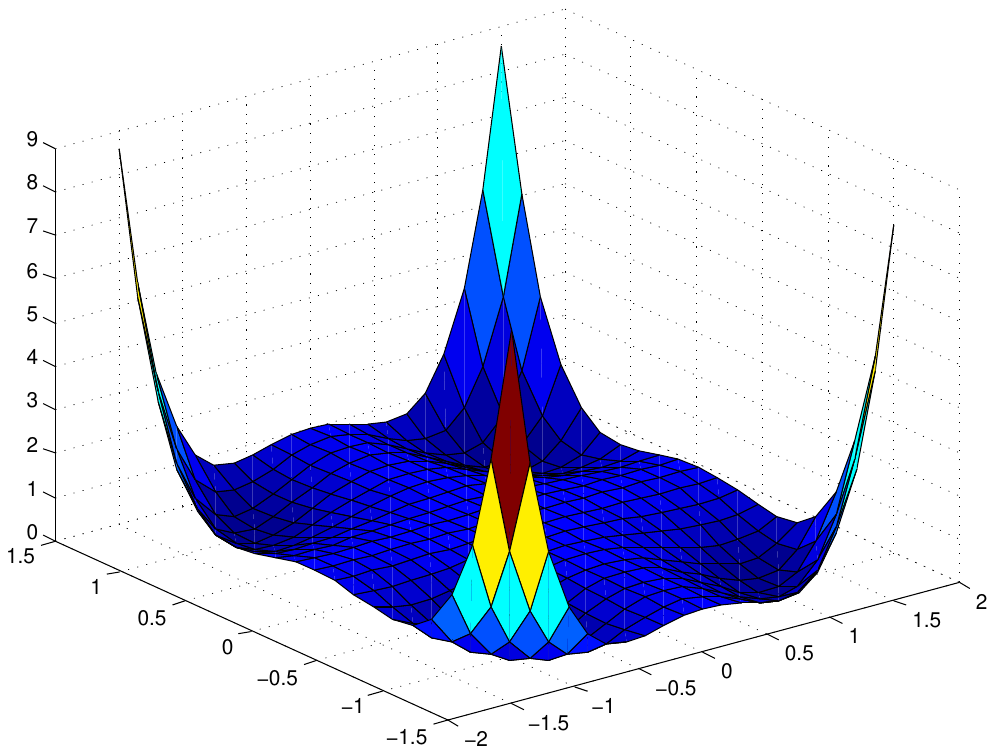}
}
\
\includegraphics[width=5cm]{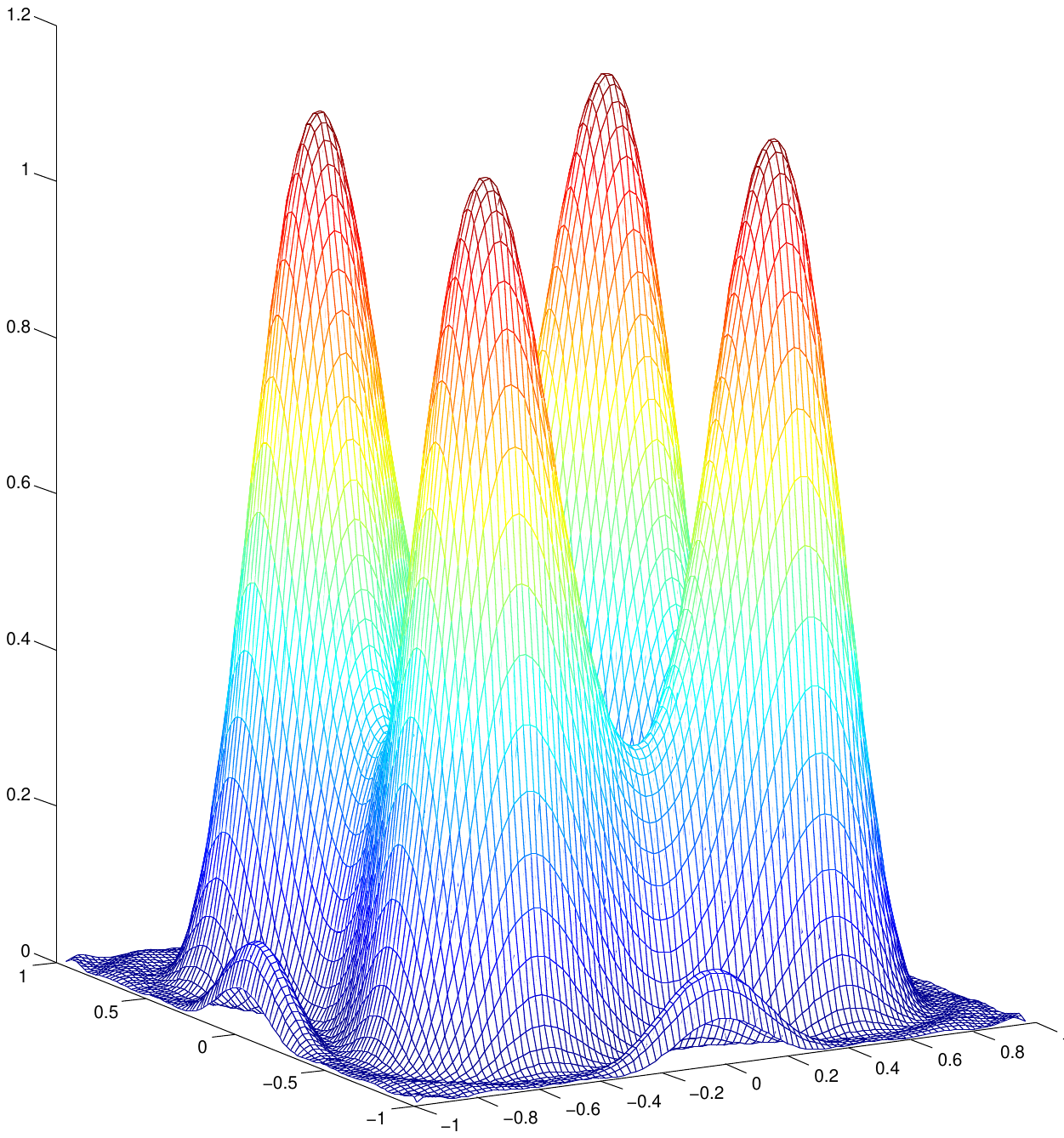}
\vspace*{-2.5cm}\caption{(a) on the left: the Motzkin polynomial (\ref{eqMotzinpol}); \ (b) on the right: the optimal degree 16 sum-of-squares polynomial $\sigma$ in program (\ref{eqfupr}).}\label{figMotzkin}
\end{figure}
\begin{trailer}{A nonnegative polynomial that is not sos}
The first explicit example of a nonnegative polynomial that is not sos was found by Motzkin in 1967, it is bivariate with degree 6 and reads
\begin{equation}\label{eqMotzinpol}
f(x,y)=x^4y^2+x^2y^4-3x^2y^2+1.
\end{equation}
It is depicted in Figure \ref{figMotzkin}(a). To see global nonnegativity one may use the arithmetic-geometric mean inequality, and use `brute force' to show $f$ is not sos (write $f$ as a sum of squares of polynomials and compare coefficients).
Also the Robinson polynomial
\begin{equation}\label{eqRobinsonpol}
f(x,y,z)= x^6+y^6+z^6- x^2y^2(x^2+y^2)-x^2z^2(x^2+z^2)-y^2z^2(y^2+z^2)+3x^2y^2z^2
\end{equation}
is globally nonnegative but not sos. We refer to Reznick \cite{Reznick2000} and Powers \cite{Powers2021} (and references therein) for a nice historic discussion and many more examples. 
\end{trailer}
At the 1900 International Congress of Mathematicians in Paris Hilbert asked whether every nonnegative polynomial can be written as a sum of squares of {\em rational} functions (known as {\em Hilbert's seventeenth problem}). This question was settled in the affirmative by Artin \cite{Artin1927}, a result that started the flourishing field of real algebraic geometry (see, e.g.,  Prestel and Delzell \cite{PrestelD2001}, Marshall \cite{Marshall2008}). 
A  landmark result is the Positivstellensatz\footnote{The terminology of {\em Positivtsellensatz} refers to  Hilbert's celebrated Nullstellenzatz that characterizes the polynomials vanishing at a given complex variety.}  by Krivine \cite{Krivine1964} and Stengle \cite{Stengle1996} that a.o. characterizes  the polynomials that are nonnegative on a   semialgebraic set $\bX$ as in (\ref{eqX}). To cite this result we need the notion of {\em preordering} $\MT(\bg)$ generated by the polynomials $\bg=\{g_1,\ldots,g_m\}$  entering the algebraic description of the set $\bX$:
\begin{equation}\label{eqTg}
\MT(\bg)=\Big\{\sum_{J\subseteq [m]} \sigma_J \prod_{j\in J}g_j: \sigma_J \in \Sigma \text{ for } J\subseteq [m]\Big\}
\end{equation}
(setting $g_\emptyset =1$).  Clearly, $\MT(\bg)\subseteq \MP(\bX)$. Krivine-Stengle show that a polynomial $f\in \oR[x]$ is nonnegative on $\bX$ if and only if $p_1 f = f^{2k}+ p_2$ for some polynomials $p_1,p_2\in \MT(\bg)$ and some integer $k\in\oN$. This leads to a sos-type decomposition of $f$ as $f=f^{2k}/p_1 + p_2/p_1$, thus `with a denominator'.

Simpler  sos-type decompositions have been shown later under more restrictive assumptions, typically assuming {\em strict positivity} of $f$. 
In particular, Reznick \cite{Reznick1995} shows a sharper result for {\em homogeneous} polynomials: If $f$ is homogeneous and satisfies $f(\bx)>0$ for all $\bx\in\oR^n\setminus \{0\}$, then there exists an integer $r\ge 0$ such that $(\sum_{i=1}^n x_i^2)^r f(\bx)\in \Sigma$.

When the set $\bX$ is compact, Schm\"udgen \cite{Schmudgen1991} shows membership in the preordering under strict positivity:
\begin{trailer}{Schm\"udgen's Positivstellensatz}
If $\bX$ is compact and 
$f$ is strictly positive on $\bX$, then $f\in \MT(\bg)$.
\end{trailer}
Under the Archimedean condition, Putinar \cite{Putinar1993} shows  membership in the {\em quadratic module} $\MQ(\bg)$, which is defined as
\begin{equation}\label{eqQg}
\MQ(\bg)=\Big\{\sum_{j=0}^m \sigma_j g_j: \sigma_j\in \Sigma \text{ for } j\in \{0,1,\ldots,m\}\Big\}
\end{equation}
(setting $g_0=1$). So, $\MQ(\bg)\subseteq \MT(\bg)\subseteq \MP(\bX)$. The advantage of the quadratic module $\MQ(\bg)$ over the preordering $\MT(\bg)$ is that it involves less terms: $m+1$ sos polynomials for $\MQ(\bg)$, instead of $2^m$ for $\MT(\bg)$.

The Archimedean condition asks that $R-\sum_{i=1}^nx_i^2\in \MQ(\bg)$ for some $R>0$ and can be seen as an algebraic certificate of compactness. 
It is a property of the algebraic description of $\bX$ rather than $\bX$ itself. It is, however,  easy to satisfy, simply  by adding the inequality of a ball containing $\bX$ to its description. The Archimedean condition holds for most sets $\bX$ of interest in applications such as the sphere, box, simplex, etc. 
\begin{trailer}{Putinar's Positivstellensatz}
If $\MQ(\bg)$ is Archimedean and 
$f$ is strictly positive on $\bX$, then $f\in \MQ(\bg)$.
\end{trailer}
So, sums of squares belong to a classical topic in real algebraic geometry, with a rich history going  back to early work by  Hilbert. It is only recently that their relevance to  optimization has been fully appreciated, starting with ground works by Shor \cite{Shor1987}, Nesterov \cite{Nesterov2000}, Lasserre \cite{Lasserre2001}, and Parrilo \cite{Parrilo2000,Parrilo2003}. A key ingredient for this link to optimization is the fact that sums of squares can be modeled using semidefinite programs, which makes them amenable to numerical algorithms.  Next, we explain how to model sos polynomials using semidefinite programming, and thereafter how to use sos polynomials to define hierarchies of (lower and upper) bounds for the original polynomial optimization problem (\ref{eqfmin}).

\subsection{Sums of Squares and Semidefinite Optimization} \label{SEC:SDP}

We begin with a quick recap on semidefinite optimization, and refer, e.g., to de~Klerk~\cite{dKbook} for a detailed treatment.  A semidefinite program   (in primal form) reads
\begin{equation}\label{eqsdpp}
p^*=\sup\{\langle C,X\rangle: \langle A_j,X\rangle =b_j \ (j\in [m]), X\in \MS^N_+\}.
\end{equation}
Here, $C,A_1,\ldots,A_m\in \MS^N$ are symmetric $N\times N$ matrices and $b\in \oR^m$ -- these are the data of the problem --  and $X\in \MS^N_+$ is the matrix variable, required to be symmetric and positive semidefinite (also written as $X\succeq 0$). 
So, the program (\ref{eqsdpp}) is a linear program over the cone $\MS^N_+$ of positive semidefinite matrices.
Its dual program reads
\begin{equation}\label{eqsdpd}
d^*=\inf\Big\{\sum_{j=1}^m b_jy_j: \sum_{j=1}^m y_jA_j-C \in \MS^N_+\Big\}.
\end{equation}
Weak duality holds: $p^*\le d^*$,  and strong duality $p^*=d^*$ holds under some Slater-type conditions. Semidefinite programs contain linear programs as a special case (when all data matrices $C,A_j$ are diagonal). 
The crucial fact is that semidefinite programs can be solved efficiently up to any given precision, under some assumptions (like knowing a small ball inside the feasible region of (\ref{eqsdpp}) and  a ball enclosing it).

As we now see, sums of squares can be modeled with semidefinite programs. As a warm-up observe that a quadratic form $f=x^T Mx$ is a sum of squares precisely when the matrix $M$ is positive semidefinite. Consider now a polynomial $f=\sum_{\alpha\in \oN^n_{2d}} f_\alpha x^\alpha $ of even degree $2d$ for which we wish to decide whether $f\in \Sigma$. 
Assume  $f\in \Sigma$, i.e.,  $f=q_1^2+\cdots +q_k^2$ for some polynomials $q_j$. Then, each $q_j\in \oR[\bx]_d$ can be written in the monomial basis as $q_j=[\bx]_d^\sfT  \ba_j$, 
where $\ba_j=((\ba_j)_\beta)_{\beta\in \oN^n_d}$ is the vector of coefficients of the polynomial $q_j$ in the monomial basis $[\bx]_d$.
 Thus, we obtain 
 $$f=\sum_{j=1}^k q_j^2= \sum_{j=1}^m [\bx]_d^\sfT \ba_j \ba_j ^\sfT [\bx]_d= [\bx]_d \Big( \sum_{j=1}^k \ba_j \ba_j^\sfT\Big) [\bx]_d
= [\bx]_d^\sfT Q [\bx]_d, $$
setting $Q:= \sum_{j=1}^k \ba_j \ba_j^\sfT$. By construction the matrix $Q$ is positive semidefinite and, by equating coefficients at both sides of the above polynomial identity, we arrive at the following characterization.
\begin{trailer}{Modeling sums of squares with semidefinite programs}
A polynomial $f=\sum_{\alpha\in\oN^n_{2d}}f_\alpha x^\alpha \in \oR[\bx]_{2d}$ is a sum of squares of polynomials if and only if the following semidefinite program has a feasible solution:
\begin{equation}\label{eqsossdp}
Q \in \MS^{N_d}_+,\quad \sum_{\beta,\gamma\in\oN^n_d: \beta+\gamma=\alpha} Q_{\beta,\gamma} = f_\alpha \ \text{ for all } \alpha\in \oN^n_{2d},
\end{equation}
where the matrix variable $Q$ is indexed by $\oN^n_d$ and $N_d=|\oN^n_d|={n+d\choose d}$.
\end{trailer}
Note that the choice of the monomial basis is not important in the above derivation, it  would work mutatis mutandis using any other polynomial basis of the polynomial space (this fact will be used in Section \ref{seceig}). 

Now that we know how to express sos polynomials using semidefinite programs we can also express polynomials in $\MQ(\bg)_{2r}$, the quadratic module truncated at a given degree $2r$, defined by
\begin{equation}\label{eqQgr}
\MQ(\bg)_{2r}=\Big\{ \sum_{j=0}^m \sigma_j g_j: \sigma_j \in \Sigma, \deg(\sigma_jg_j)\le 2r \text{ for } j=0,1,\ldots,m\Big\},
\end{equation}
or in the truncated preordering $\MT(\bg)_{2r}$ defined by
\begin{equation}\label{eqpreordering2r}
\MT(\bg)_{2r}=\Big\{\sum_{J\subseteq [m]} \sigma_J \prod_{j\i J}g_j: \sigma_J\in\Sigma, \deg\Big(\sigma_J\prod_{j\in J}g_j\Big)\le 2r \text{ for } J\subseteq [m]\Big\}.
\end{equation}
For this, set $d_j=\lceil \deg(g_j)/2\rceil$, so that polynomials in the truncated quadratic module $\MQ(\bg)_{2r}$ can be modeled as 
$\sum_{j=0}^m g_j [\bx]_{r-d_j}^\sfT Q_j [\bx]_{r-d_j}$ for some positive semidefinite matrices $Q_j\in \MS^{N_{d_j}}$ (for $j=0,,\ldots,m$).
In other words, membership in the truncated quadratic module can be modeled as a semidefinite program. The same holds of course for membership in the truncated preordering.

We now have all tools in hands to define the hierarchies of upper bounds and lower bounds for the polynomial optimization problem (\ref{eqfmin}). 

\subsection{Upper Bounds} 

To define  upper bounds on $\fmin$  the starting point is to observe that problem (\ref{eqfmin}) can be reformulated as a linear optimization problem over the set $\MM(\bX)$ of positive Borel measures supported within the set $\bX$:
\begin{equation}\label{eqfminmeas}
\fmin =\min\Big\{\int_\bX f(\bx)d\nu(\bx): \nu\in \MM(\bX),\ \int_\bX d\nu(\bx)=1\Big\}.
\end{equation}
The argument is simple: On the one hand, $\int_\bX f d\nu\ge \fmin$ since $\nu$ is a probability measure on $\bX$. On the other hand, the Dirac delta  $\nu=\delta_{\bx^*}$  at a global minimizer $\bx^*$ of $f$ in $\bX$ provides a feasible solution to (\ref{eqfminmeas}) with value $f(\bx^*)=\fmin$.

Let is now fix  a (reference) measure $\mu\in \MM(\bX)$ whose support is equal to $\bX$. The next step is  that one can restrict the optimization in (\ref{eqfminmeas}) to the measures $\nu$ that have a sos density with respect to this given measure $\mu$:
\begin{equation}\label{eqfmindensity}
\fmin=\min\Big\{ \int_\bX f(\bx)\sigma(\bx)d\mu(\bx): \sigma\in \Sigma,\ \int_\bX \sigma(\bx)d\mu(\bx)=1\Big\}
\end{equation}
(Lasserre \cite{Lasserre2010}). Intuitively, this relies  on the fact that the Dirac delta $\delta_{\bx^*}$  can be well approximated  by sos polynomials.
The upper bounds on $\fmin$ are then obtained by restricting the optimization to sos polynomials of a given degree.

\begin{trailer}{The upper bounds on $\fmin$}
For any $r\in \oN$, define the parameter
\begin{equation}\label{eqfupr}
\ub(f,\bX,\mu)_r=\min \Big\{ \int_\bX f(\bx)\sigma(\bx)d\mu(\bx): \sigma\in \Sigma_{2r},\ \int_\bX \sigma(\bx)d\mu(\bx)=1\Big\}.
\end{equation}
We have $\fmin \le \ub(f,\bX,\mu)_{r+1}\le \ub(f,\bX,\mu)_r$.
\end{trailer}
As an illustration, we show in Figure \ref{figMotzkin}(a) the Motzkin polynomial and in (b) the optimal sos density $\sigma$ for program (\ref{eqfupr}) at order $2r=16$ (which approximates well the sum of Dirac functions at the four minimizers $(\pm 1, \pm 1)$).

So, the parameter (\ref{eqfupr}) depends on the choice of the reference measure $\mu$ on $\bX$ and  it  can be computed via a semidefinite program (in fact, as an eigenvalue problem as we see in Section \ref{seceig}). In view of (\ref{eqfmindensity}), the bounds $\ub(f,\bX,\mu)_r$ converge asymptotically to $\fmin$ as $r\to\infty$. We return to these bounds in Section \ref{secupperbound}, where we will discuss in detail the convergence rate of  the error range $\ub(f,\bX,\mu)_r-\fmin$  as the relaxation order $r$ grows. 

For now, we summarize in Table \ref{TAB:ub} the known results about this convergence rate for various classes of compact sets $X$ (and reference measures $\mu$).
 As we see there, one can show a convergence rate in $O(1/r^2)$ (up to a $\log^2(r)$ factor) 
 for a large class of compact sets. 
 For instance, one can show a rate in $O(\log^2(r)/r^2)$ for any convex body $\bX$, and in $O(1/r^2)$ when $\bX$ is `round', which means that there are tangent inner and outer balls at any point on its boundary. 
 The convergence rate in $O(1/r^2)$  is essentially optimal, as can be seen for the case $\bX=[-1,1]$ (due to the explicit link to extremal roots of orthogonal polynomials, see relation (\ref{eqxub}) in Section \ref{secroots}), or for the case $\bX=\mathbb S^{n-1}$ (due to a link to cubature rules, see Section \ref{seccubature}).
 
  \begin{table}
{\renewcommand{\arraystretch}{1.2}%
\begin{small}
\begin{tabular}{cccc}
$\mainset$ \textbf{(compact)} & \textbf{error}  & \textbf{measure} $\mu$ & \textbf{reference}\\
\hline 
\makecell{Geometric assumption} & $O(1/\sqrt{r})$  & Lebesgue & de Klerk, Laurent, Sun \cite{deKlerkLaurentSun2017} \\
\makecell{Convex body} & $O(1 / r)$  & Lebesgue & de Klerk, Laurent \cite{deKlerkLaurent2018} \\
\makecell{Semialgebraic} & $O(\log^2(r)/r^2)$  & Lebesgue & Slot, Laurent \cite{SlotLaurent2021b} \\
\makecell{with dense interior,} &   &  & \\
\makecell{convex body} &   &  & \\
\hline
$\mathbb S^{n-1}$ 	  & $O(1/r)$  & uniform & Doherty, Wehner \cite{DohertyWehner2013} \\
$\mathbb S^{n-1}$ 	  & $O(1/r^2)$  & uniform & de Klerk, Laurent \cite{deKlerkLaurent2022} \\
$[-1, 1]^n$   & $O(1/r^2)$  & \makecell{$\prod_{i}(1-\bx_i)^\lambda d\x~(\lambda = -\frac{1}{2})$} &de Klerk, Laurent  \cite{deKlerkLaurent2020} \\
$[-1, 1]^n$   & $O(1/r^2)$  & \makecell{$\prod_{i}(1-\bx_i)^\lambda d\x~(\lambda \geq -\frac{1}{2})$} & Slot, Laurent \cite{SlotLaurent2022a} \\
`Round' convex body & $O(1/r^2)$ & Lebesgue & Slot, Laurent \cite{SlotLaurent2022a} \\
$\ball{n}$ 	  & $O(1/r^2)$  & $(1-\|\bx\|^2)^\lambda d\x~(\lambda \geq 0)$ & Slot, Laurent\cite{SlotLaurent2022a} \\
$\simplex{n}$ & $O(1/r^2)$  & Lebesgue & Slot, Laurent \cite{SlotLaurent2022a}  \\
 \\
\hline
\end{tabular}
\end{small}
}
\caption{Overview of known results on the asymptotic error $\ub(f,\bX,\mu)_r-\fmin$ of Lasserre's 
 hierarchy of upper bounds. }
\label{TAB:ub}
\end{table}  
 
\subsection{Lower Bounds}

To define lower bounds on $\fmin$ the starting point is to observe that the minimum value taken by the polynomial $f$ over the set $\bX$ is equal to the largest scalar $\lambda$ for which the polynomial $f-\lambda$ is nonnegative over $\bX$:
\begin{equation}\label{fminpos}
\fmin=\sup\{\lambda: \lambda\in \oR,\ f(\bx)-\lambda\ge 0 \ \text{ for all } \bx\in \bX\}.
\end{equation}
Now, we obtain bounds by replacing the positivity condition by membership in either the truncated quadratic module $\MQ(\bg)_{2r}$ or preordering $\MT(\bg)_{2r}$ (see \eqref{eqQgr}, \eqref{eqpreordering2r}).
\begin{trailer}{The lower bounds on $\fmin$}
For any integer $r\ge \deg(f)/2$, define the parameters
\begin{align}
\label{eqflowr}
\lb(f,\MQ(\bg))_r &= \sup\{\lambda: \lambda\in \oR,\ f-\lambda \in \MQ(\bg)_{2r}\}, \\
\label{eqflowrschmud}
\lb(f,\MT(\bg))_r &= \sup\{\lambda: \lambda\in \oR,\ f-\lambda \in \MT(\bg)_{2r}\}.
\end{align}
By definition, we have $\lb(f,\MQ(\bg))_r \leq \lb(f,\MT(\bg))_r$ for all $r$. Furthermore, we have $\lb(f,\MQ(\bg))_r\le \lb(f,\MQ(\bg))_{r+1} \le \fmin$, and $\lb(f,\MT(\bg))_r\le \lb(f,\MT(\bg))_{r+1} \le \fmin$. 
\end{trailer}
As an application of the earlier mentioned Positivstellens\"atze of Putinar and Schm\"udgen,  if $\MQ(\bg)$ is Archimedean (resp., $\bX$ is compact), then the bounds $\lb(f,\MQ(\bg))_r$ (resp., $\lb(f,\MT(\bg))_r$) converge to $\fmin$ as $r\to \infty$ (Lasserre \cite{Lasserre2001}).
We return to these bounds in Section \ref{seclowerbound}, where we will discuss their asymptotic convergence rates as the relaxation order $r$ grows. We summarize the main known results in Table \ref{TAB:lb} below. As we see, the results depend on the algebraic structure of the semialgebraic set $\bX$.

\begin{table}
\centering
{\renewcommand{\arraystretch}{1.2}%
\begin{tabular}{cccc}
$\mainset$ \textbf{(compact)} & \textbf{error} & \textbf{certificate} & \textbf{reference}\\
\hline
Archimedean   & $O(1/\log(r)^c)$ & $\qmod{\bg}$ & Nie, Schweighofer \cite{NieSchweighofer} \\
{Archimedean}   & {$O(1/r^c)$} & {$\qmod{\bg}$} & {Baldi, Mourrain \cite{BaldiMourrain2023}, Baldi et al. \cite{BaldiMourrainParusinski2022}} \\
General  	  & $O(1/r^c)$ & $\preo{\bg}$ & Schweighofer \cite{Schweighofer2004} \\
$[-1, 1]^n$   & $O(1/r)$ & $\qmod{\bg}$ & Baldi, Slot \cite{BaldiSlot2024} \\
\hline
$\mathbb S^{n-1}$ 	  & $O(1/r^2)^\star$ & $\qmod{\bg} ~(= \preo{\bg})$ %
& Fang, Fawzi  \cite{FangFawzi2021} \\
$\ball{n}$ 	  & $O(1/r^2)$ & $\qmod{\bg} ~(=\preo{\bg})$ & Slot \cite{Slot2022} \\
$\simplex{n}$ & $O(1/r^2)$ & $\preo{\bg}$ & Slot \cite{Slot2022} \\
$[-1, 1]^n$   & $O(1/r^2)$ & $\preo{\bg}$ & Laurent, Slot  \cite{LaurentSlot2023} \\
\hline
\end{tabular}
}
\caption{Overview of known results on the asymptotic error $\fmin -\lb(f,\bX)_r$ of Lasserre's hierarchies of lower bounds. ($^\star$ The result on the sphere $\oS^{n-1}$ is shown in \cite{FangFawzi2021} only for homogeneous polynomials~$f$, {but it can be extended to general polynomials \cite{BLS}}.)}
\label{TAB:lb}
\end{table}

\subsection{Link to Cubature Rules} \label{seccubature}
A  natural approach to design upper bounds on the minimum value of $f$ over $\bX$ is by minimizing $f$ over a well-chosen finite set of points $\bX_r\subseteq \bX$. 
For example, for  a set $\bX \subseteq [0,1]^n$, one may take for $\bX_r$ the set of rational points in $\bX$ that have a given  denominator $r\ge 1$. This clearly gives a hierarchy of upper bounds that converge asymptotically to $\fmin$ when the denominator $r$  grows. For the simplex $\bX=\Delta^n$, 
 we have $|\bX_r|={n+r\choose r}$, thus polynomial in $n$ for any fixed $r$. In fact, the parameters $f_{\text{\rm min},\bX_r}$   lead to a polynomial-time approximation scheme with convergence rate 
$\min_{\bx\in \bX_r}f(\bx)- \fmin= O(1/r)$ (de Klerk, Laurent, Parrilo \cite{dKLP2006}). 
However, for the hypercube $\bX=[0,1]^n$, we have $|\bX_r|= (r+1)^n$, which is thus {\em exponential} in the dimension $n$. We refer to Martinez et al. \cite{MPSV2019} for related work for well-chosen  finite meshes on $\bX$. Let us mention the following simple link to cubature rules.

Let $\mu\in\MM(\bX)$  with support $\bX$. Assume that $\bX_r=\{x^{(i)}: i\in [N]\}\subseteq \bX$, together with positive weights $w=(w_i)_{i\in [N]}$, provides a positive cubature rule for $(\bX,\mu)$ that is exact at degree $\deg(f)+2r$. That is, for any polynomial $p$ with $\deg(p)\le \deg(f)+2r$ we have 
$\int_\bX pd\mu=\sum_{i=1}^N w_i p(x^{(i)})$.
Then, one  easily sees that 
$$\ub(f,\bX,\mu)_r\ge \min_{\bx \in \bX_r} f(\bx)= f_{\text{\rm min}, \bX_r}  \ge \fminX.
$$
Sometimes, this permits to show tightness of the convergence rate for the upper bounds $\ub(f,\bX,\mu)_r$ from information about cubature rules. This is done, e.g.,  for the sphere $\bX=\mathbb S^{n-1}$ (equipped with the uniform Haar measure) by de Klerk, Laurent \cite{deKlerkLaurent2022}, who show $\ub(f,\bX,\mu)_r-\fmin=\Omega(1/r^2)$ for linear $f$.

\section{Performance Analysis of the Upper Bounds}\label{secupperbound}

Recall the set up: $\bX\subseteq \oR^n$ is a compact set and $\mu$ is a reference measure supported by $\bX$. Consider the bilinear form   $\langle \cdot,\cdot\rangle_\mu$ on $\oR[\bx]$ induced by $\mu$, defined by
\begin{equation}\label{eqinnerprodmu}
\langle p,q\rangle _\mu=\int_\bX p(\bx)q(\bx)d\mu(\bx)\ \text{ for } p,q\in \oR[\bx].
\end{equation}
If $\bX$ has a nonempty interior, then  this defines an inner product on $\oR[\bx]$. Otherwise, if $\bX$ has an empty interior (e.g., $\bX=\oS^{n-1}$ is the unit sphere), then this defines an inner product on  the space $\MR(\bX)$ of polynomials on $\bX$. Let $\MB_\mu=\{P_\alpha: \alpha \in \oN^n\}$ be an orthonormal basis of $\MR(\bX)$ with respect to this inner product, with the property that the set 
\begin{equation}\label{eqbasismu}
\MB_{\mu,d}=\{P_\alpha: \alpha\in \oN^n_d\}, \ \ \text{sometimes also denoted as the vector  } [P_\alpha]_d,
\end{equation}
 is a basis of $\MR(\bX)_d$ (the set of polynomials that agree on $\bX$ with a degree $d$ polynomial)  for each $d\in \oN$.

Our objective here is to give a (rough) overview of the main ideas used to show the results in Table \ref{TAB:ub}, where rates in $O(1/r^2)$ and in $O(\log^2(r)/r^2)$ are presented.

The starting point for showing a  $O(r^2)$ rate is establishing an eigenvalue reformulation for the parameter $\ub(f,\bX,\mu)_r$ (Section \ref{seceig}) and a link to extremal roots of orthogonal polynomials in the univariate case (Section \ref{secroots}). In addition, one uses some simple `tricks' to reduce the analysis to (at most) quadratic polynomials and to simpler sets and measures (Section \ref{sectricks}).

To show the (slightly weaker) rate $O(\log^2 (r)/r^2)$ for much more general sets $\bX$ one  again follows a univariate strategy, now reducing the search in problem (\ref{eqfupr}) to sos polynomials $\sigma(\bx)=s(f(\bx))$, where $s$ is a {\em univariate} sos polynomial (see Section \ref{secunilog}).

\subsection{Some Useful  Tricks for the Analysis}\label{sectricks}
We group here some simple facts, useful for the analysis of the upper bounds. 

\begin{lemma}[de Klerk, Laurent, Sun \cite{deKlerkLaurentSun2017}] \label{lemestimator}
Let $f, g\in \oR[\bx]$ satisfy (i) $f(\bx)\le g(\bx)$ for all $\bx \in \bX$, and (ii)  $f$ and $g$ take the same minimum value on $\bX$, i.e., $\fminX=\gminX$. Then, for any $r\in \oN$, we have
$\ub(f,\bX,\mu)_r-\fminX\le \ub(g,\bX,\mu)_r-\gminX.
$
\end{lemma}
Using Taylor's expansion, one can see that any polynomial $f$ admits an upper estimator $g$  on $\bX$ that satisfies the assumptions (i), (ii) of Lemma \ref{lemestimator} and that is (linear or) quadratic. 

The next lemma permits to reduce the analysis for a pair $(\bX,w)$, where $\bX$ is equipped with the absolutely continuous measure $w(\bx)d\bx$,  to a possibly simpler pair $(\bhX,\hw)$, where $\bX\subseteq \bhX$ and $w\le \hw$ `look the same' around a minimizer.
 
\begin{lemma}[Slot, Laurent \cite{SlotLaurent2022a}]\label{lemsimpleset}
Let $\bX\subseteq \bhX$ be compact sets, where $\bX$ (resp., $\bhX$) is equipped with an absolutely continuous measure $w(\bx)d\bx$ (resp., $\hw(\bx)d\bx$). 
Let $x^*\in \bX$ be a minimizer of $f$ in $\bX$. Assume the following conditions hold: 
\begin{description}
\item[(i)] $\bX$, $\bhX$ are `locally similar' at $\bx^*$: $\bX\cap B^n(\bx^*,\epsilon) = \bhX\cap B^n(\bx^*,\epsilon)$ for some $\epsilon>0$ (with $B^n(\bx^*,\epsilon) $ the ball centered at $\bx^*$ with radius $\epsilon$).
\item[(ii)] $w(\bx)\le \hw(\bx)$ for all $\bx\in \text{int}(\bX)$.
\item[(iii)]  $w$, $\hw$ are `locally comparable' at $\bx^*$: $C\cdot \hw(\bx)\le w(\bx)$ for $\bx\in \text{int}(K)\cap B^n(\bx^*,\epsilon)$ for some $\epsilon>0$ and $C>0$.
\end{description} 
Then, there exists a (linear or) quadratic polynomial $g$ such that $f(\bx)\le g(\bx)$ for all  $\bx\in \bhX$, $\fminX=\gminX$, and $\ub(f,\bX,w)_r -\fminX \le {2\over C}(\ub(g,\bhX,\hw)_r-\gminhX)$.
\end{lemma}

The up-shot of these two lemmas  can be summarized as follows.
\begin{trailer} {\bf Recipe for analyzing the upper bounds} It suffices to analyze the convergence rate  of the range  $\ub(f,\bX,\mu)_r-\fmin $ for {\em (linear or) quadratic}  $f$ and `simple' sets $\bX$ (like the box $[-1,1]^n$ or the ball $B^n$).
\end{trailer}

\subsection{Reformulation as an Eigenvalue Problem}\label{seceig}

We begin with giving an eigenvalue reformulation for the parameter $\ub(f,\bX,\mu)_r$ from (\ref{eqfupr}). For this, we express the  polynomial $\sigma\in \Sigma_{2r}$ entering the definition of the parameter $\ub(f,\bX,\mu)_r$ using
 the orthonormal basis $\MB_{\mu,r}$, ordered as the vector $[P_\alpha]_r$ (as in \eqref{eqbasismu}). So, we can write 
$\sigma=\langle Q, [P_\alpha]_r[P_\alpha]_r^\sfT\rangle$, where 
$Q$ is the matrix variable (indexed by $\MB_{\mu,r}$). Then, we have 
$\int_\bX f\sigma d\mu= \langle Q, M_{\mu,r}(f)\rangle$ and $\int_\bX \sigma d\mu=\Tr(Q)$, after defining the matrix $$M_{\mu,r}(f)=\Big(\int_\bX f P_\alpha P_\beta d\mu\Big)_{\alpha,\beta\in \oN^n_r}.$$
In this way we arrive at the following {\em eigenvalue reformulation}:
\begin{equation}\label{equpeig}
\ub(f,\bX,\mu)_r= \min\{ \langle M_{\mu,r}(f), Q\rangle: Q\succeq 0,\ \Tr(Q)=1\}= \lambda_{\min}(M_{\mu,r}(f)),
\end{equation}
showing that  the parameter $\ub(f,\bX,\mu)_r$ is equal to the smallest eigenvalue of the matrix $M_{\mu,r}(f)$.
Some remarks are in order here: this computation relies on the matrix $M_{\mu,r}(f)$, which in turn requires to be able to integrate a polynomial on the set $\bX$ w.r.t. the measure $\mu$.
Thus, for practical computation  one needs to restrict to some relatively easy sets equipped with well-understood measures. 

Estimating this smallest eigenvalue remains a difficult problem in general. However, there is a situation where it is very well-understood: in the univariate case for the interval $\bX=[-1,1]$ equipped with a `nice' measure. We first consider  this case, which will form the basis for understanding  the general multivariate case.

\subsection{Univariate Case: Links to Roots of Orthogonal Polynomials}\label{secroots}

Here, we consider the univariate case $n=1$ and the interval $\bX=[-1,1]$ equipped with a measure $\mu$ supported on $[-1,1]$. Consider as above the orthonormal basis $\MB_\mu=\{P_k: k\ge 0\}$ of $\oR[x]$ w.r.t. the inner product $\langle \cdot,\cdot\rangle _\mu$.  Then, the polynomials $P_k$ satisfy the well-known {\em 3-term recurrence}: there exist scalars $a_k,b_k$ ($k\ge 0$) such that
\begin{equation}\label{eq3term}
xP_k= a_{k-1}P_{k-1} +b_kP_k+ a_k P_{k+1} \ \text{ for } k\ge 0 \quad (\text{setting } a_{-1}=0).
\end{equation}
Observe that for the polynomial $f=x$,  the matrix $M_{\mu,r}(x)=\big(\int_{-1}^1xP_iP_jd\mu\big)_{i,j=0}^r$ (also known as the  {\em Jacobi matrix}) is tri-diagonal and  its  eigenvalues are the roots of the degree $r+1$ orthogonal polynomial $P_{r+1}$. In particular,
\begin{equation}\label{equproot}
\lambda_{\min}(M_{\mu,r}(x)) = \text{ smallest root of the  orthogonal polynomial } P_{r+1} \text{ (w.r.t.  } \mu).
\end{equation}
We refer, e.g., to Dunkl and Xu \cite{DunklXu2014} for background on orthogonal polynomials.

From now on, we will assume that the reference measure $\mu $ is of {\em Jacobi-type}.
\begin{trailer}{Jacobi-type measures on $[-1,1]$}%
The Jacobi-type measure is $d\mu(x)= (1-x)^\lambda (1+x)^{\lambda'}dx$, where   $\lambda,\lambda'>-1$; the associated orthogonal polynomials are known as the Jacobi polynomials.\\
The case when  $\lambda=\lambda'=-1/2$ is of special interest; then, the 3-term recurrence (\ref{eq3term}) reads $xP_k= (P_{k-1}+P_{k+1})/2$ and the  associated orthogonal polynomials are the Chebychev polynomials.  \\
The reason for restricting to  Jacobi-type  measures is that the behaviour of the extremal roots of their associated orthogonal polynomials is well-understood:  The smallest root of $P_{r}$ is $-1 +\Theta(1/r^2)$ (Dimitrov, Nikolov \cite{DN2010}, Driver, Jordaan \cite{DJ2012}).
\end{trailer}
As observed earlier, it suffices  to analyze the upper bounds for {\em linear} or {\em quadratic}~$f$. 

For {\em linear} $f=x$, since  the smallest root of $P_{r+1}$  is $-1+\Theta( 1/r^2)$, by  combining with relations  (\ref{equpeig})  and (\ref{equproot}), we directly obtain  the  error analysis 
\begin{equation}\label{eqxub}
\ub(x,[-1,1],\mu)_r -\fmin= \lambda_{\min}(M_{\mu,r}(x))+1 =\Theta\Big({1\over r^2}\Big)
\end{equation}
(from de Klerk, Laurent \cite{deKlerkLaurent2020}).

We now consider the case when $f$ is {\em quadratic}. We distinguish two cases, depending whether $f$ has a minimizer in the boundary or in the interior of $[-1,1]$. In the former case, $f$ admits a linear upper estimator; then  we  can apply Lemma \ref{lemestimator} and use  (\ref{eqxub}) to conclude that $\ub(f,[-1,1],\mu)_r -\fmin =O(1/r^2)$. 

So, assume now $f=x^2 +c x$ with $a\in (-2,2)$, so that $f$ attains its minimum value  in $(-1,1)$. Then, one needs to estimate the smallest eigenvalue of the matrix $M_{\mu, r}(f)$, a difficult task in general. Indeed, $M_{\mu, r}(f)$ is  a 5-diagonal matrix, with its entries depending on the parameters $a_k,b_k$ in the 3-term recurrence (\ref{eq3term}) and the parameter $c$ in the definition of $f$.  Estimating this smallest eigenvalue is, however, easier in the Chebyshev case.

In the case when the measure $\mu$ is of Chebyshev-type, i.e.,  $d\mu=(1-x^2)^{-1/2}dx$, the  matrix $M_{\mu,r}(f)$ is `almost' a circulant matrix: After modifying its first two rows and columns it can be made a circulant matrix, whose eigenvalues can be explicitly computed. Combining with an interlacing argument, one can show again that $\ub(f,[-1,1],\mu)_r -\fmin =O(1/r^2)$ 
 (de Klerk and Laurent \cite{deKlerkLaurent2020}).
Combining this with using Lemma \ref{lemsimpleset} one can show the same result for a general measure $\mu$ with weight $(1-x^2)^\lambda$ when $\lambda\ge -1/2$ (Slot and Laurent \cite{SlotLaurent2022a}). 

In summary, we have $\ub(f,[-1,1],\mu)_r-\fmin  =O(1/r^2)$ for any polynomial $f$, in the case when $[-1,1]$ is equipped with a measure $\mu$ with Jacobi weight $(1-x^2)^\lambda$ and $\lambda\ge -1/2$.
This error estimate in $O(1/r^2)$ extends then to the box $[-1,1]^n$, equipped with a product of such measures. 

Similarly, to establish the rate $O(1/r^2)$ for the other sets in Table \ref{TAB:ub}, one  uses the  result (\ref{eqxub}) for the univariate case, combined with the results of Lemmas \ref{lemestimator} and \ref{lemsimpleset} (see Slot, Laurent \cite{SlotLaurent2022a} for the ball, simplex, and round convex bodies), and possibly some `integration trick' (see de Klerk, Laurent \cite{deKlerkLaurent2022} for the sphere).

\subsection{Another Analysis Technique using Needle Polynomials}\label{secunilog}

At this point there remains to explain how to show the (slightly weaker) rate $O(\log^2(r)/r^2)$ presented in Table \ref{TAB:ub}  for  general compact sets $\bX$ like semialgebraic sets (with dense interior) and convex bodies. For this the key idea is again to follow a univariate approach, as explained by Slot, Laurent \cite{SlotLaurent2021b}. Instead of searching over all (multivariate) sos polynomials $\sigma\in \oR[\bx]$ in program (\ref{eqfupr}), the idea is to restrict the search to {\em univariate} sos polynomials $s \in \oR[x]$ and then set $\sigma(\bx)=s(f(\bx))$. In this way one arrives at the following (weaker) bounds (introduced by Lasserre \cite{Lasserre2021}):
\begin{align}\label{eqfubpf}
\begin{split}
\ub_{\#}(f,\bX,\mu)_r & =\min\Big\{ \int_\bX f(\bx) s(f(\bx)) d\mu(\bx): s\in \oR[x]_{2r}  \text{ sos},\ \int_\bX s(f(\bx))d\mu(\bx)=1\Big\}\\
& = \min\Big\{\int_{f(\bX)} x s(x) df_{\#}\mu(x): s\in \oR[x]_{2r} \text{ sos},\ \int_{f(\bX)}s(x) df_{\#}\mu(x)=1\Big\},
\end{split}
\end{align}
where $f_{\#}\mu$ denotes the univariate measure obtained by taking the push-forward of $\mu$ by $f$. Then, we have $\ub(f,\bX,\mu)_{rd}\le \ub_{\#}(f,\bX,\mu)_r$ if $d=\deg(f)$.

In view of relations (\ref{equpeig}) and (\ref{equproot}), the analysis of the bounds 
$\ub_{\#}(f,\bX,\mu)_r$ relies on the smallest roots of the orthogonal polynomials w.r.t. the push-forward measure $f_{\#}\mu$. However, these orthogonal polynomials are not well understood in general, so another approach is needed.
Up to translation  one may assume that $\fmin=0$. Then, the idea is to find a univariate sos polynomial $s$ that approximates well the Dirac delta at the origin. For this, one can employ the so-called {\em needle polynomials} (from Kro\'o \cite{Kroo}) that are widely used in the literature of approximation theory. We refer to Slot, Laurent \cite{SlotLaurent2021b} for the technical details and extension to compact sets satisfying a suitable geometric assumption.

\section{Performance Analysis of the Lower Bounds} \label{seclowerbound}
We turn now to the lower bounds. Let $\bX = \{\bx\in \oR^n: g_1(\bx)\ge 0,\ldots, g_m(\bx)\ge 0\}$ be as in~\eqref{eqX}. We assume that the quadratic module $\MQ(\bg)$ is Archimedean. Recall that this implies that, for any $f$, the bounds $\lb(f, \MQ(\bg))_r$ and $\lb(f, \MT(\bg))_r$ converge to $\fmin$ as $r \to \infty$ by Putinar's and Schm\"udgen's Positivstellens\"atze, respectively. In this section, we outline the main ideas used to obtain the convergence rates presented in~Table~\ref{TAB:lb}. For ease of writing, we will use the letter $\MC$ to refer to either $\qmod{\bg}$ or $\preo{\bg}$. A useful observation is the following.
\begin{observation} \label{OBS:representation}
Let $f \in \oR[\bx]$. Then, for any $\epsilon > 0$ and $r \in \oN$, we have
\[
	\fmin - \lb(f, \MC)_r \leq \epsilon \iff f - \fmin + \epsilon \in \MC_{2r}.
\]
\end{observation}
Thus, in order to prove convergence rates for Lasserre's hierarchies of lower bounds, it suffices to find low-degree sos-representations of $f - \fmin + \epsilon$. In what follows, we discuss two methods of obtaining such representations. We give particular attention to the so-called \emph{polynomial kernel method} in Section~\ref{SEC:PKM} below, as it reveals an interesting connection between the analysis of the \emph{lower} bounds and the \emph{upper} bounds discussed in Section~\ref{secupperbound}.

Throughout this section, we set $d=\deg(f)$, which should be thought of as being fixed, while the relaxation order $r$ will grow.

\subsection{Analysis via the Polynomial Kernel Method} \label{SEC:PKM}
The \emph{polynomial kernel method} (PKM) has been successful for proving strong convergence rates for optimization over certain distinguished sets $\bX$, including $\bX = \oS^{n-1}, B^n, \Delta^n, [-1, 1]^n$.
It was initially described\footnote{Earlier, weaker analyses of the lower bounds on $\oS^{n-1}$ due to Reznick~\cite{Reznick1995} and Doherty, Wehner~\cite{DohertyWehner2013} already relied on the PKM \emph{implicitly}.} in the special case $\bX = \oS^{n-1}$ by Fang, Fawzi~\cite{FangFawzi2021}, who use it to show a convergence rate in $O(1/r^2)$ in that setting. 
It was later applied by Laurent, Slot~\cite{LaurentSlot2023} and Slot~\cite{Slot2022}, to prove rates in $O(1/r^2)$ for optimization over $[-1, 1]^n$ and  $\Delta^n, B^n$, respectively. The latter work by Slot~\cite{Slot2022} is the first to describe the technique in full generality, and we follow its exposition here.

Suppose we were able to construct a linear operator $\kernelop : {\oR[\bx] \to \oR[\bx]}$ with the following three properties:
\begin{align}
\tag{P1} \label{EQ:P1}
\kernelop(1) &= 1, \\
\tag{P2} \label{EQ:P2}
\kernelop p &\in \MC_{2r} \quad \text{for all } p \in \MP(\bX), \\ 
\tag{P3} \label{EQ:P3}
\max_{\bx \in \bX} | \kernelop^{-1} f(\bx) - f(\bx) | &\leq \epsilon. 
\end{align}
Then, we claim that $(f - \fmin) + \epsilon \in \MC_{2r}$ (which shows that $\fmin - \lb(f, \MC)_r \leq \epsilon$ by Observation~\ref{OBS:representation}). 
Indeed, by \eqref{EQ:P3}, we have $\kernelop^{-1}f(\bx)\ge f(\bx)-\epsilon\ge \fmin-\epsilon$ on $\bX$. Using \eqref{EQ:P1}, we obtain $\kernelop^{-1} (f-\fmin+\epsilon)= \kernelop^{-1}f -\fmin+\epsilon\in \MP(\bX)$.
By~\eqref{EQ:P2}, we may then conclude that 
\begin{equation*}
	f - \fmin + \epsilon = \kernelop \big [\kernelop^{-1} \big(f - \fmin + \epsilon\big) \big] \in \MC_{2r}.
\end{equation*}

\subsubsection{Constructing Linear Operators}
It remains to construct operators $\kernelop$ that enjoy these special properties. For this, we rely on the theory of (polynomial) \emph{reproducing kernels}. Let $\mu$ be a (sufficiently nice) measure supported on $\bX$. Then, any polynomial $\kernel \in \oR[\bx, \by]$ induces a linear operator $\kernelop$ on $\oR[\bx]$ via convolution:
\begin{equation} \label{EQ:convolution}
	\kernelop p (\bx) := \int_{\bX} \kernel(\bx, \by) p(\by) d \mu(\by) \quad (p \in \oR[\bx]).
\end{equation}
The goal is to choose $\kernel$ in such a way that $\kernelop$ satisfies \eqref{EQ:P1}-\eqref{EQ:P3}. For \eqref{EQ:P2}, it suffices to choose $\kernel$ so that the polynomial $\bx \mapsto \kernel(\bx, \by)$ lies in $\MC_{2r}$ for all \emph{fixed} $\by \in \bX$. This is a consequence of  Tchakaloff's Theorem, which allows us to write the integral~\eqref{EQ:convolution} as (positively weighted) finite sum over a cubature rule (see~Laurent, Slot \cite{LaurentSlot2023}). 
To establish \eqref{EQ:P1} and~\eqref{EQ:P3}, it turns out that it is enough to control the \emph{eigenvalues} of $\kernelop$. 
The eigenvalues of $\kernelop$ can be related to the polynomial $\kernel$ directly if we make the assumption that it is of the form
\begin{equation} \label{EQ:expansion}
	\kernel(\bx, \by) = \sum_{\|\alpha\|_1 \leq 2r} \lambda_{\alpha} P_\alpha(\bx) P_\alpha(\by) \quad (\lambda_\alpha \in \oR),
\end{equation}
where $\{ P_\alpha : \alpha \in \oN^n \}$ is an orthonormal basis of $\oR[\x]$ w.r.t. $\mu$, ordered so that $\mathrm{deg}(P_\alpha) = \|\alpha\|_1$ for all $\alpha$. Indeed, by orthonormality, the eigenvalues of $\kernelop$ are then given by the coefficients $\lambda_\alpha$ in~\eqref{EQ:expansion}. If we set $\lambda_\alpha = 1$ for all $\alpha$, the resulting polynomial is called the \emph{reproducing kernel} for $\oR[\bX]_{2r}$ w.r.t. $\mu$. Its associated operator $\kernelop$ acts as the identity on $\oR[\bx]_{2r}$. It thus satisfies \eqref{EQ:P1} and \eqref{EQ:P3} (for $\epsilon = 0$). However, it does not have an sos-representation in general, and so it does not satisfy~\eqref{EQ:P2}. The idea is to \emph{perturb} the reproducing kernel slightly, choosing $\lambda_\alpha \approx 1$ (for all $|\alpha|\le d$), in such a way that all three properties hold. The following lemma makes the relation between the $\lambda_\alpha$ and properties~\eqref{EQ:P1}, \eqref{EQ:P3} precise.
\begin{lemma} \label{LEM:gamma}
	Let $\kernel$ be as in~\eqref{EQ:expansion}, with $\lambda_0 = 1$ and $\lambda_\alpha \in (1/2, 1]$ for all $\|\alpha\|_1 \leq d$. Then, its associated linear operator $\kernelop$ satisfies~\eqref{EQ:P1}, and
	\[
		\max_{\bx \in \bX} | \kernelop^{-1} f(\bx) - f(\bx) | \leq \gamma \cdot \sum_{|\alpha| \leq d} (1-\lambda_\alpha) \cdot (f_{\max} - f_{\min}).
	\]
	for some constant $\gamma > 0$ depending only on $\bX$ and $d$ (but not on $f$).
\end{lemma}
The inequality in Lemma~\ref{LEM:gamma} above can be improved in certain special cases (such as $\bX = \mathbb S^{n-1}$). In particular, the constant $\gamma$ can be bounded nontrivially, but this is beyond the scope of this survey.

\subsubsection{A Connection to the Upper Bounds}

 In general, it is not obvious how to choose the $\lambda_\alpha$ such that \eqref{EQ:P2} holds. However, for the aforementioned special choices of $\bX$, one can rely on special structure of the reproducing kernel to make this problem tractable.
For instance, in the original application of the PKM by Fang, Fawzi \cite{FangFawzi2021} for the sphere  $\oS^{n-1}$, the classical \emph{Funk-Hecke formula} tells us that, for any $k \in \oN$,
\[
\sum_{\|\alpha\|_1 = k} P_\alpha(\bx) P_\alpha(\by) = \mathcal{G}_k(\bx \cdot \by) \quad (\bx, \by \in \oS^{n-1}),
\]
where $\mathcal{G}_k$ is the \emph{Gegenbauer polynomial} of degree $k$. These are the (univariate) orthogonal polynomials on the interval $[-1, 1]$ w.r.t to weight $w(t) = (1-t^2)^{\frac{n-3}{2}}$. Assuming that the $\lambda_\alpha$ in~\eqref{EQ:expansion} depend only on $|\alpha|$, the polynomial $\kernel$ then equals
\[
	\kernel(\bx, \by) = \sum_{k=0}^{2r} \lambda_k \mathcal{G}_k(\bx \cdot \by).
\]
Write $q(t) := \sum_{k=0}^{2r} \lambda_k \mathcal{G}_k(t)$. Note that if $q \in \Sigma[t]_{2r}$ is a (univariate) sum of squares, then $\kernel(\bx, \by)$ lies in the quadratic module of the sphere for fixed $\by \in \oS^{n-1}$, and so $\kernelop$ will satisfy \eqref{EQ:P2}. 
In light of Lemma~\ref{LEM:gamma}, we thus wish to select the $\lambda_k$ so that $\lambda_0 =1$~\eqref{EQ:P1}; $q \in \Sigma[t]_{2r}$~\eqref{EQ:P2}; and $\sum_{k=1}^d (1 - \lambda_k)$ is as small as possible~\eqref{EQ:P3}. Remarkably, the optimal selection of the $\lambda_k$ corresponds to a particular instance of Lasserre's \emph{upper} bounds. This correspondence allows us to transport the analysis of these bounds in Section~\ref{secupperbound} to our present setting. Indeed, by orthogonality, and after choosing the right normalization of the Gegenbauer polynomials~$\mathcal{G}_k$, we have that 
$
{\lambda_k = \int_{-1}^1 \mathcal{G}_k(t) q(t) w(t)dt}.
$  
Now set $g(t) = d - \sum_{k=0}^d \mathcal{G}_k(t)$. Then, we have that 
\[
\int_{-1}^1 g(t) q(t) w(t) dt = \sum_{k=1}^d (1 - \lambda_k).
\]
Thus, choosing the $\lambda_k$ optimally reduces to solving the optimization problem
\[
	\mathrm{opt} := \inf_{q \in \Sigma[t]_{2r}} \bigg \{ \underbrace{\int_{-1}^1 \sum_{k=1}^{d} g(t) q(t) w(t) dt}_{\sum_{k=1}^d (1 - \lambda_k)} : \underbrace{\int_{-1}^1 q(t) w(t) dt}_{\lambda_0} = 1 \bigg \}.
\]
We recognize this as the program~\eqref{eqfupr} that defines the upper bound $\ub({g, [-1,1]})_r$ for the minimization of $g$ on $[-1,1]$. As $g_{\min} = 0$ (attained at $t=1$), we may conclude that $\mathrm{opt} = O(1/r^2)$, and a convergence rate of the \emph{lower} bounds on $\oS^{n-1}$ of the same order follows. We refer to Fang, Fawzi~\cite{FangFawzi2021} (and the exposition in Slot \cite{Slotthesis}) for details.

For $\bX \in \{B^n, \Delta^n\}$, one has (more complicated) analogs of the Funk-Hecke formula. This allows one to establish a correspondence between upper and lower bounds similar to the above.  For $\bX = [-1,1]^n$, there is no such formula. However, the reproducing kernel has a product structure in that case, which can be exploited to obtain representations of $\kernel$ in the truncated preordering. We refer to Slot~\cite{Slot2022, Slotthesis} and Laurent, Slot~\cite{LaurentSlot2023} for details.

\begin{trailer}{Recipe for analyzing the lower bounds (I): Polynomial kernel method}%
Convergence rates for Lasserre's hierarchy of lower bounds may be obtained from sum-of-squares representations of (perturbed) reproducing kernels on~$\bX$. For distinguished semialgebraic sets $\bX$, such representations can be found via the hierarchy of \emph{upper} bounds.
\end{trailer}

\subsection{Analysis via Algebro-Geometric Reduction}
We now discuss a technique employed first by Schweighofer~\cite{Schweighofer2004}, Nie, Schweighofer \cite{NieSchweighofer} and later by Baldi, Mourrain~\cite{BaldiMourrain2023} and Baldi et al.~\cite{BaldiMourrainParusinski2022}, to prove convergence rates for Lasserre's lower bounds on general (compact, Archimedean) semialgebraic sets $\bX = \Se(\bg)$ (as in \eqref{eqX}). With respect to the polynomial kernel method, this technique yields weaker guarantees, but it does not rely on special structure of the set $\bX$. We follow the exposition of Baldi, Mourrain~\cite{BaldiMourrain2023}.

Suppose that, for a given polynomial $f$, we wish to show that $f - \fmin + \epsilon \in \qmod{\bg}_{2r}$ for some small $\epsilon > 0$. The idea is to embed the (potentially complicated) set $\bX$ in a simple semialgebraic set $\bY = \Se(\bh) \supseteq \Se(\bg) = \bX$. For instance, $\bY$ might be a (scaled) ball, simplex or box. Following Baldi, Mourrain~\cite{BaldiMourrain2023}, let us fix $\Se(\bh) = [-1, 1]^n$ for concreteness. Using specialized results on the convergence of Lasserre's hiearchy on $[-1, 1]^n$ (e.g., those obtained via the PKM), we know that, for some $\epsilon = O(1/r^2)$,
\begin{equation} \label{EQ:specializedh}
	f - f_{\min, \bY} + \epsilon \in \preo{\bh}_{2r}.
\end{equation}
This fact is not immediately useful to us, for two reasons. First, we wish to obtain a representation in $\qmod{\bg}$, not $\preo{\bh}$. Second, the minimum $f_{\min, \bY}$ of $f$ on $\bY$ might be (much) smaller than $\fmin = f_{\min, \bX}$. 
The first issue is resolved by proving an inclusion $\preo{\bh}_{2r} \subseteq \qmod{\bg}_{2r + \ell}$. Note that the existence of such an inclusion is not surprising, as $\bY \supseteq \bX$, and thus $\MP(\bX) \supseteq \MP(\bY)$. Nonetheless, careful arguments are needed to control the additional degree term $\ell$, but we do not discuss these here. The second issue is more serious. It is resolved as follows:
Suppose we were able to construct a regularizing polynomial $q \in \qmod{\bg}_{\ell'}$ with the property that 
$f_{\min, \bX} \leq (f - q)_{\min, \bY}$. Then, setting $\tilde f = f - q$, we could apply~\eqref{EQ:specializedh} to~$\tilde f$, and use the inclusion of 
$\preo{\bh}$ in $\qmod{\bg}$  to obtain, for some $\tilde \epsilon = O(1/r^2)$,
\[
	\tilde f - \tilde f_{\min, \bY} + \tilde \epsilon \in \preo{\bh}_{2r} \subseteq \qmod{\bg}_{2r + \ell}.
\]
Using the identity
$f-f_{\min,\bX} +\tilde \epsilon =( \tilde f -\tilde f_{\min,\bY} +\tilde \epsilon ) + q + \tilde f_{\min,\bY}-f_{\min,\bX}$, combined with $\tilde f_{\min, \bY} \geq f_{\min, \bX}$ and $q \in \qmod{\bg}_{\ell'}$,
this gives us  the representation 
\[
	f - f_{\min, \bX} + \tilde \epsilon \in \qmod{\bg}_{\max\{2r + \ell, \, \ell'\}}.
\]
Note that $\tilde \epsilon$ depends on the degree of $q$, which may be substantially larger than $d$. Therefore, $\tilde \epsilon \gg \epsilon$ in general. In addition to the extra degree terms $\ell, \ell'$, this is what explains the weaker rates obtained from this technique compared to the PKM.

\subsubsection{Construction of the Regularizing Polynomial}

It remains to construct the polynomial $q \in \qmod{\bg}$ described above. As mentioned, this construction has a large impact on the quality of the resulting convergence guarantees. In fact, it is arguably the construction of $q$ that sets apart the results of Nie, Schweighofer~\cite{NieSchweighofer} and Baldi et al.  \cite{BaldiMourrain2023, BaldiMourrainParusinski2022}, leading to stronger rates in the latter works. We give a sketch of the construction in Baldi, Mourrain~\cite{BaldiMourrain2023}.

Recall that after setting $\tilde f = f - q$, we want that $\tilde f_{\min, \bY} \geq f_{\min, \bX}$. Thus, we wish that $q \ll 0$ on $\bY \setminus \bX$, and $q \approx 0$ on $\bX$ (as $q \in \qmod{\bg}$, we know that $q \geq 0$ on $\bX$). 
Without loss of generality, assume that the polynomials $g_j$ defining $\bX$ satisfy $-1 \leq g_j(\x) \leq 1$ for all $\x \in \bY$ and $j \in [m]$.
The idea is to first construct a \emph{univariate} sum-of-squares polynomial $h \in \Sigma[t]$, with the property that, for some small $\eta > 0$, we have $h(t) \approx 0$ for $t \in [\eta, 1]$ and $h(x) \approx 1$ for $t \in [-1, -\eta]$. This can be done via Chebyshev approximation, see Baldi, Mourrain~\cite[Section 2.3]{BaldiMourrain2023}. Then, we set $q(\bx) = M \cdot \sum_{j=1}^m g_j(\bx) h(g_j(\bx))$, with $M > 0$. Note that $q \in \qmod{\bg}$ by definition. Note further that, for $\bx \not \in \bX$ sufficiently far from the boundary $\mathrm{bd}(\bX)$, at least one $g_j(\bx) < -\eta$, whence $g_i(\bx) h(g_i(\bx)) \ll 0$, whereas for $\x \in \bX$ sufficiently far from $\mathrm{bd}(\bX)$, all $g_j(\bX) \in [\eta, 1]$, whence $q(\bx) \approx 0$. It remains to analyze the situation for $\bx$ close to $\mathrm{bd}(\bX)$, where (some of the) constraints $g_j(\bx)$ are close to 0. There, a careful comparison of the gradients $\nabla f$ of $f$ and $\nabla g_j$ of the constraints $g_j$ is required to finish the argument. This comparison involves the so-called \L{ojasiewicz}-constant of~$\bg$, which appears in the exponent of the final convergence rate. See~Baldi, Mourrain~\cite{BaldiMourrain2023} for details.

\begin{trailer}{Recipe for analyzing the lower bounds (II): Algebro-geometric reduction}%
Convergence rates for Lasserre's hierarchy of lower bounds on general semialgebraic sets may be derived from rates on special sets (such as the simplex or hypercube), through delicate algebraic and geometric arguments.
\end{trailer}

\subsection{Extensions}
We end this section with some recent extensions of the techniques described above.

\subsubsection{Algebro-geometric reduction for the hypercube}
Baldi and Slot~\cite{BaldiSlot2024} use a variation of the reduction technique described above to prove convergence rates in $O(1/r)$ for the Putinar-type lower bounds on $[-1, 1]^n$, based on the quadratic module (recall that Laurent, Slot~\cite{LaurentSlot2023} only show rates for the stronger Schm\"udgen-type bounds, based on the preordering). This extension shows that the algebro-geometric method can also be useful for analyses on specific, rather than general sets $\bX$. The main difference w.r.t. the above is that the `simple' set $\bY \supseteq [-1, 1]^n$ used in Baldi, Slot~\cite{BaldiSlot2024} depends on $f, r$ (while it was static before). In fact, somewhat remarkably, $\bY = \Se(\bh) = [1-\eta, 1+\eta]^n$ is itself just a scaled hypercube, with $\eta = \eta(f, r) > 0$. The use of a dynamic $\bY$ eliminates the need for the regularizing polynomial~$q$, permitting to prove stronger rates than those in Baldi et al.~\cite{BaldiMourrain2023, BaldiMourrainParusinski2022} (which would give a rate in $O(1/\sqrt[10]{r})$ in this setting). On the other hand, proving an inclusion $\preo{\bh} \subseteq \qmod{\bg}$ with proper degree bounds is more complicated. See Baldi, Slot~\cite{BaldiSlot2024} for details.

\subsubsection{Sparse polynomial optimization}
As we have seen in Section~\ref{SEC:SDP}, the semidefinite programs used to model Lasserre's hierarchies in $n$ variables at degree $r$ involve matrices of size $N_r = {n + r\choose r}$. These SDPs are thus intractable already for moderately large values of $n$. To address this issue, variants of Lasserre's hierarchy of lower bounds that exploit \emph{sparsity} of the underlying polynomial optimization problem have been proposed in the literature, see Magron et al.~\cite{Magron:sparsebook} for an overview. A recent work by Korda et al.~\cite{KordaMagronRZ2024} extends both of the techniques discussed in Section~\ref{seclowerbound} to this setting, yielding a performance analysis for the sparse bounds. Furthermore, the convergence rates achieved there are stronger than those of Section~\ref{seclowerbound} for sufficiently sparse POPs (relative to the size of the resulting SDPs).

\subsubsection{Generalized moment problems}
Polynomial optimization is a special case of the \emph{generalized problem of moments} (GPM), which asks to minimize a linear function over the cone of positive Borel measures on $\oR^n$ under some linear conditions on the moments. Lasserre's approximation hierarchies naturally extend to the GPM, see Lasserre \cite{Lasserre2009}, de Klerk, Laurent~\cite{deKlerkLaurent:survey}. Convergence rates for these extended hierarchies have been established, especially in the context of dynamical systems, optimal control and volume estimation, by Korda et al.~\cite{KordaHenrionJones2017, KordaHenrion2018} and Schlosser et al.~\cite{STL2024, ST2024}. 
The analyses in the latter two recent works 
rely directly on (a combination of) existing convergence results for polynomial optimization discussed above.

\section{Discussion}

\subsubsection*{Tightness of the performance analysis}
As we have seen, the convergence rates for the hierarchy of \emph{upper} bounds presented in Section~\ref{secupperbound} are essentially tight (up to log-factors). For the \emph{lower} bounds, the situation is much less clear. The literature on `worst-case' examples for the lower bounds is mostly \emph{qualitative} in nature, see, e.g., Scheiderer~\cite{Scheiderer2005}, Powers, Scheiderer~\cite{PowersScheiderer2001}. 
On the quantitative side, Baldi and Slot~\cite{BaldiSlot2024} recently showed that the Putinar-type lower bounds on $\bX = [-1, 1]^n$ converge at a rate no better than $\Omega(1/r^8)$. Note that this negative result is still rather far away from the best-known positive result in that setting, which is in $O(1/r)$. It is an interesting research direction to prove negative results on the convergence of the lower bounds that either match the positive results more closely; or apply to more general $\bX$; or apply to the Schm\"udgen-type bounds.

\subsubsection*{Exponential convergence under local optimality conditions}

A common feature of the convergence rates in Table~\ref{TAB:lb} is that they are \emph{sub}exponential, in particular no better than $O(1/r^c)$ for some constant $c>0$ depending on the set~$\bX$ and the type of certificate. In light of the negative results mentioned above, these are likely `optimal' (up to improving the exponent $c$). On the other hand, Bach and Rudi~\cite{BachRudi} show an \emph{exponential} convergence rate for the Schm\"udgen-type bounds on~$ \bX = [-1, 1]^n$, which holds under an additional assumption on the objective function $f$. Roughly speaking, $f$ should have a (strictly) positive definite Hessian at its global minimizer $\bx^* \in [-1, 1]^n$, see Bach, Rudi~\cite{BachRudi} for details.\footnote{Note that this assumption holds {\em generically}. 
Under a similar assumption, Nie~\cite{Nie2014} showed earlier that the hierarchy of lower bounds has finite convergence for general $\bX=\Se(\bg)$, i.e., $f-\fmin$ belongs to $\MQ(\bg)_r$ for some $r \in \oN$.  However, his result gives no quantitative information on $r$.}
Under this assumption, $f - \fmin$ can be written as a sum of squares of \emph{smooth} functions (not necessarily polynomials). In turn, these smooth functions can be approximated by polynomials via (truncated) expansion in the Fourier basis, leading to a sum-of-squares representation of $f - \fmin$. The exponential convergence rate for Lasserre's lower bounds then follows from the fact that the Fourier coefficients of a smooth function decay exponentially quickly. 
As the approach of Bach and Rudi relies primarily on Fourier analysis, it seem likely that it could be extended to other distinguished $\bX$, such as $\bX = \oS^{n-1}$. More unclear is whether it could also be used to prove guarantees for general semialgebraic sets~$\bX$, which we believe is an interesting open question.

\subsubsection*{Optimization over finite semialgebraic sets}
In this survey, we have focused on the setting where $\bX \subseteq \oR^n$ is an infinite set. Sum-of-squares hierarchies have been extensively studied for optimization over finite semialgebraic $\bX$ as well, particularly for (subsets of) the boolean hypercube: $\bX \subseteq \{-1, 1\}^n$ (or $\bX\subseteq \{0,1\}^n$). A key difference in the finite setting is that the hierarchies always converge in a finite number of steps (under a minor condition on the description of~$\bX$), see Nie \cite{Nie2013}. In fact, they converge in  $\Theta(n)$ steps for a semialgebraic subset of the boolean cube, see Laurent \cite{Laurent2003}, Fawzi et al. \cite{FSP2016}.
For this reason, asymptotic analysis as $r \to \infty$ does not make sense. Rather, one often fixes the level $r \in \oN$ of the hierarchy, and lets the number of variables $n$ tend to infinity. There is a large literature on this regime in the theoretical computer science community, see, e.g., Barak, Steurer~\cite{BarakSteurerSurvey} and references therein. One could also consider a `hybrid' regime, where $r = c \cdot n$ for some constant $c > 0$ and $n, r \to \infty$ simultaneously. In this regime, the polynomial kernel method (see Section~\ref{SEC:PKM}) yields error guarantees for the upper and lower bounds on $\bX = \{-1, 1\}^n$, see Slot, Laurent~\cite{SlotLaurent2022b} for details.

\end{document}